# A proof of Lindelöf's hypothesis
## by
## M. Aslam Chaudhry

Department of Mathematics and Statistics
King Fahd University of Petroleum and Minerals
Dhahran 31261, Saudi Arabia.
E-mail: maslam@kfupm.edu.sa

## Abstract
We define an associated Lindelöf-function for the ratio of the zeta functions and use its representation to get a unique extension of Lindelöf's function that proves Lindelöf's hypothesis.

## 1. Introduction

Lindelöf's hypothesis remains one of the most outstanding open problems in mathematics. It is concerned with the growth of the Riemann zeta function in the critical strip $\{s := \sigma + i\tau : 0 \leq \sigma \leq 1, -\infty < \tau < \infty\}$. Riemann proved that the zeta-function

$$\zeta(s) := \sum_{n=1}^{\infty} \frac{1}{n^s} \qquad (s = \sigma + i\tau, \sigma > 1), \qquad (1.1)$$

has a meromorphic continuation to the complex plane. This satisfies the functional equation (see [7], p.13 (2.1.1))

$$\chi(s) := \zeta(s)/\zeta(1-s) = \frac{1}{\pi}(2\pi)^s \sin(\frac{\pi}{2}s)\Gamma(1-s), \qquad (1.2)$$

and has simple zeros at $s = -2, -4, -6, \ldots$ called the *trivial zeros*. All the other zeros, called the *non-trivial zeros*, of the function are symmetric about the *critical line* $\sigma = 1/2$ in the *critical strip* $0 \leq \sigma \leq 1$. The multiplicity of these non-trivial zeros (in general) is not known. Riemann conjectured that the non-trivial zeros of the function lie on the critical line $\sigma = 1/2$. This conjecture is called the Riemann hypothesis. Many numerical experiments support the hypothesis. However, an analytic proof of the hypothesis is still needed. There are several equivalent formulations of the hypothesis ([7], pp.282-328) that depend on the estimate of $|\zeta(\sigma + i\tau)|$ ($|\tau| \to \infty$). It is to be remarked that the function $\chi(s)$ ($0 < \sigma < 1$) is non-zero, analytic and $\chi(s)\chi(1-s) = 1$. Moreover, we have the inverse Mellin transform representation ([6], p. 91(3.3.6))

$$2\cos(2\pi x) = \frac{1}{2\pi i} \int_{c-i\infty}^{c+i\infty} (2(2\pi)^{-s} \cos(\frac{\pi s}{2})\Gamma(s))x^{-s}ds \qquad (0 < c < 1/2). \qquad (1.3)$$

The line of integration, as explained in ([6], p. 91), can not be moved outside the half strip $0 < \sigma < 1/2$. By using the Riemann functional equation (1.2), we can rewrite (1.3) as

$$\lambda(x) := 2\cos(2\pi x) = \frac{1}{2\pi i} \int_{c-i\infty}^{c+i\infty} \chi(1-s) x^{-s} ds \qquad (0 < c < 1/2), \qquad (1.4)$$

which is a useful representation. In particular, it follows from the inverse Mellin transform relations (1.4) that $\chi(1-s)$ is an absolutely and uniformly bounded function in every sub-strip $(0 < \sigma_1 \leq \sigma \leq \sigma_2 < 1/2)$. Since $\chi(s)\chi(1-s) = 1$ and $\chi(s)$ $(0 < \sigma < 1)$ is non-zero, therefore, it follows from (1.4) that $\chi(s)$ is absolutely and uniformly bounded in every sub-strip $(1/2 < \sigma_1 \leq \sigma \leq \sigma_2 < 1)$ and

$$\frac{1}{x}\lambda(\frac{1}{x}) := \frac{2}{x}\cos(2\pi/x) = \frac{1}{2\pi i} \int_{c-i\infty}^{c+i\infty} \chi(s) x^{-s} ds \qquad (1/2 < c < 1), \qquad (1.5)$$

which is a useful representation for the study of the properties of the zeta function. We have introduced the associated Lindelöf-function for the $\chi$-function and find its representation using (1.4) and (1.5). The function is found to be non-negative, continuous, non-increasing and convex downwards and vanishes at $\sigma = 1/2$. It was shown by Lindelöf, by using the modified theorem ([3], p. 186), that

$$|\zeta(\sigma + i\tau)| < K\tau^{\frac{1}{2}(1-\sigma)} \log \tau \qquad (0 \leq \sigma \leq 1, \tau > 1). \qquad (1.6)$$

Lindelöf's hypothesis says ([7], p. 276) that for every $\varepsilon > 0$,

$$|\zeta(\sigma + i\tau)| = O(\tau^\varepsilon) \qquad (\sigma \geq 1/2, |\tau| \to \infty). \qquad (1.7)$$

The above hypothesis was formulated by E. Lindelöf in 1908. Let $\mu(\sigma)$ be the least upper bound of the number $A$ such that $|\zeta(\sigma + i\tau)|\tau^{-A}$ is bounded when $\tau \to \infty$ ([3], p. 186). The Lindelöf hypothesis is equivalent to the hypothesis that

$$\mu(\sigma) = 0 \qquad (\sigma \geq 1/2). \qquad (1.8)$$

Since the function $\mu(\sigma)$ is non-negative, continuous, non-increasing and convex downwards, the relation (1.8) is equivalent to the statement that $\mu(1/2) = 0$. Huxley [5] has shown that $\mu(1/2) \leq 0.1561$ and this seems to be the best known result in this direction. Backlund showed in 1918 ([3], p.188) that Lindelöf's hypothesis follows from the Riemann hypothesis. There are several other equivalent formulations of the Lindelöf hypothesis (see [4], [5] and [7]). One of the most important equivalent formulations of the hypothesis is ([7], p. 276),

$$\frac{1}{T}\int_1^T \left|\zeta(\frac{1}{2} + i\tau)\right|^{2k} d\tau = O(T^\varepsilon) \qquad (\forall \varepsilon > 0, T \to \infty, k = 1,2,3,...). \qquad (1.9)$$

The result (1.9) was proven for $k=1$ by Littlewood and for $k=2$ by Heath-Brown [5]. For $k=3$, the result is still open. The function $\mu(\sigma)$ satisfies the functional equation ([3], p. 186)

$$\mu(\sigma) - \mu(1-\sigma) = \frac{1}{2} - \sigma \qquad (-\infty < \sigma < \infty). \tag{1.10}$$

It is to be remarked that the functional equation is also satisfied by the functions $\mu_A(\sigma) := \frac{1}{2}(1-\sigma)$ and $\mu_B(\sigma) := \frac{1}{2}(1-\sigma)^2$ ($-\infty < \sigma < \infty$). Thus the functional equation does not have a unique solution. We define an associated Lindelöf-function for the ratio of the zeta functions and show that it does satisfy Lindelöf's hypothesis. We use its representation and exploit the Riemann functional equation to get a unique extension of Lindelöf's function that proves Lindelöf's hypothesis. Our method of proof is extended to prove the generalized Lindelöf's hypothesis.

## 2. Proof of Lindelöf's hypothesis

**Theorem** The Lindelöf hypothesis is true.

**Proof** Putting $s = i\tau$ ($\tau \geq \tau_0 > 0$) in (1.7), we find

$$\chi(i\tau) = \frac{\tau^2}{2}(2\pi)^{i\tau}[\sinh(\pi\tau/2)/(\pi\tau/2)]\Gamma(-i\tau) \tag{2.1}$$

which leads to the representation

$$|\chi(i\tau)| = \frac{\tau^2}{2}[\sinh(\pi\tau/2)/(\pi\tau/2)]|\Gamma(-i\tau)|$$
$$= \frac{\tau}{\pi}[\sinh(\pi\tau/2)]|\Gamma(-i\tau)|. \tag{2.2}$$

However, we have ([1], p. 4(1.19)) the well known representation

$$|\Gamma(i\tau)| = \frac{\sqrt{\pi}}{\sqrt{\tau\sinh(\pi\tau)}} \qquad (\tau \geq \tau_0 > 0). \tag{2.3}$$

From (2.2) and (2.3), we find

$$|\chi(i\tau)| = \frac{\sqrt{\tau}}{\sqrt{\pi}}[\sinh(\pi\tau/2)]/\sqrt{\sinh(\pi\tau)} = \frac{\sqrt{\tau}}{\sqrt{2\pi}}\left(\frac{1-e^{-\pi\tau}}{\sqrt{1-e^{-2\pi\tau}}}\right)$$
$$(\tau \geq \tau_0 > 0). \tag{2.4}$$

However, as

$$\frac{1-e^{-x}}{(1-e^{-2x})^{1/2}} = \frac{e^x - 1}{(e^{2x}-1)^{1/2}} = \sqrt{\frac{e^x - 1}{e^x + 1}} \leq \sqrt{\frac{e^x}{e^x + 1}} \leq 1 \quad (x \geq x_0 > 0) \text{ is bounded above by one}$$

and it approaches one as $x \to \infty$, the equation (2.4) leads to the sharp inequality

$$|\chi(i\tau)| \leq \frac{\sqrt{\tau}}{\sqrt{2\pi}} \qquad (\tau \geq \tau_0 > 0), \tag{2.5}$$

that holds for $\tau \geq 0$ as well. From the identity $\chi(s)\chi(1-s)=1$, we find

$$\left|\chi(\frac{1}{2}+i\tau)\right|=1 \qquad (\tau \geq 0). \tag{2.6}$$

We have shown that the analytic function $\chi(\sigma+i\tau)$ ($0 \leq \sigma \leq 1/2, \tau \geq \tau_0 > 0$) satisfies the relations (2.5) and (2.6). Therefore, according to modified Lindelöf's theorem ([3], p. 186), we have a closed form of the corresponding affine function $k(\sigma)$ ($p=1/2$, $q=0$, $\sigma_1=0$, $\sigma_2=1/2$)

$$k(\sigma) = [(q-p)/(\sigma_2 - \sigma_1)](\sigma - \sigma_1) + p = 1/2 - \sigma, \tag{2.7}$$

such that

$$\left|\chi(\sigma+i\tau)\right| \leq K\tau^{k(\sigma)} = K\tau^{\frac{1}{2}-\sigma} \qquad (0 \leq \sigma \leq 1/2, \tau \geq \tau_0 > 0), \tag{2.8}$$

which is sharp and consistent with the known asymptotic relation ([7], p. 81)

$$\left|\chi(\sigma+i\tau)\right| \sim (\frac{\tau}{2\pi})^{\frac{1}{2}-\sigma} \qquad (\tau \to \infty), \tag{2.9}$$

in any fixed strip $\alpha \leq \sigma \leq \beta$. Moreover, we have shown in Appendix-A that in (2.8), the absolute constant $K \leq 8$ for all $\tau \geq 1$.

Let $\mu_\chi(\sigma)$ ($-\infty < \sigma < \infty$) be Lindelöf's function associated with the $\chi$-function. We define $\mu_\chi(\sigma)$ as the least upper bound of the non-negative numbers $A$ such that $\left|\chi(\sigma+i\tau)\right|\tau^{-A}$ is bounded when $\tau \to \infty$. It is to be remarked that the requirement of the number $A$ to be non-negative in the above definition is not needed for the zeta function as it is always non-negative due to the known properties (see [3], pp.184-186) of the zeta function. The relation (2.9) shows that the function $\mu_\chi(\sigma)$ ($-\infty < a \leq \sigma \leq b < \infty$) is well defined, non-negative, continuous and bounded. We note from (1.10) and (1.11) that $\left|\chi(\sigma+i\tau)\right|$ ($1/2 < \sigma_1 \leq \sigma \leq \sigma_2 < 1, \tau \geq \tau_0 > 0$) is absolutely and uniformly bounded. Moreover, $\left|\chi(1/2+i\tau)\right|=1$ is on the critical line. Therefore, from the above analysis and by using (2.8) and (2.9), we find (see [7], p. 81(5.1.2))

$$\mu_\chi(\sigma) = 0 \qquad (1/2 \leq \sigma \leq \sigma_2 < \infty), \tag{2.10}$$

$$\mu_\chi(\sigma) = \frac{1}{2} - \sigma \qquad (-\infty < \sigma \leq 1/2). \tag{2.11}$$

Hence the function $\mu_\chi(\sigma)$ ($-\infty < \sigma < \infty$) satisfies Lindelöf's hypothesis. Let $H(\sigma-a)$ be the unit step function defined by

$$H(\sigma-a) = \begin{cases} 1 & \sigma > a, \\ 0 & \sigma < a. \end{cases} \tag{2.12}$$

We note that

$$H(\sigma-a) + H(a-\sigma) = 1 \qquad (\sigma \neq a). \tag{2.13}$$

The step function $H(\sigma - a)$ is discontinuous at $\sigma = a$ where it can be defined arbitrarily. We define $H(0) := c$ where $c \in (0, 1/2)$ for our purposes which leads to $H(\sigma - a) + H(a - \sigma) \leq 1$ ($-\infty < \sigma < \infty$). We can express the associated Lindelöf function in (2.10) and (2.11) uniquely as

$$\mu_\chi(\sigma) = (\frac{1}{2} - \sigma) H(\frac{1}{2} - \sigma) \qquad (-\infty < \sigma < \infty), \qquad (2.14)$$

which shows that $\mu_\chi(\sigma)$ is continuous, non-increasing, convex downwards and non-negative. We note that $\mu_\chi(1/2) = 0 H(0) = 0(c) = 0$ is independent of $c \in (0, 1/2)$. Moreover, it satisfies both Lindelöf's hypothesis and the functional equation as

$$\mu_\chi(\sigma) - \mu_\chi(1 - \sigma) = (\frac{1}{2} - \sigma)[H(\frac{1}{2} - \sigma) + H(\sigma - \frac{1}{2}) = \frac{1}{2} - \sigma$$
$$(-\infty < \sigma < \infty). \qquad (2.15)$$

Let us write Lindelöf's function $\mu_\zeta(\sigma) := \mu(\sigma)$ that is associated with the $\zeta$-function. From the Riemann functional equation

$$\zeta(s) = \chi(s)\zeta(1 - s) \qquad (-\infty < \sigma < \infty), \qquad (2.16)$$

and due to the fact that $\zeta(s)$ ($\sigma \geq 1 + \delta$, $\forall \delta > 0$) is bounded (see [7], p. 81 (5.1.2)), we find that

$$\mu_\zeta(\sigma) = \mu_\chi(\sigma) = 1/2 - \sigma \qquad (-\infty < \sigma \leq 0), \qquad (2.17)$$

and

$$\mu_\zeta(\sigma) \leq \mu_\chi(\sigma) + \mu_\zeta(1 - \sigma) \qquad (-\infty < \sigma < \infty). \qquad (2.18)$$

From (2.14), (2.17) and (2.18), we get

$$\mu_\zeta(\sigma) = \mu_\chi(\sigma) = (\frac{1}{2} - \sigma) H(\frac{1}{2} - \sigma) \qquad (-\infty < \sigma \leq 0). \qquad (2.19)$$

and

$$\mu_\zeta(\sigma) \leq (\frac{1}{2} - \sigma) H(\frac{1}{2} - \sigma) + \mu_\zeta(1 - \sigma) \qquad (-\infty < \sigma < \infty). \qquad (2.20)$$

It is to be remarked that the presence of the factor $H(\frac{1}{2} - \sigma)$ in (2.19) and (2.20) is important. Using (2.13) we can rewrite (2.20) as ($\sigma \neq 1/2$),

$$\mu_\zeta(\sigma) \leq (\frac{1}{2} - \sigma + \mu_\zeta(1 - \sigma)) H(\frac{1}{2} - \sigma) + \mu_\zeta(1 - \sigma) H(\sigma - \frac{1}{2})$$
$$(-\infty < \sigma < \infty), \qquad (2.21)$$

From (1.10) and (2.21), we find that ($\sigma \neq 1/2$),

$$\mu_\zeta(\sigma) \leq \mu_\zeta(\sigma) H(\frac{1}{2} - \sigma) + \mu_\zeta(1 - \sigma) H(\sigma - 1/2)$$
$$(-\infty < \sigma < \infty), \qquad (2.22)$$

which can be rewritten as ($\sigma \neq 1/2$),

$$\mu_\zeta(\sigma)(1 - H(1/2 - \sigma)) \leq \mu_\zeta(1 - \sigma) H(\sigma - 1/2) \qquad (-\infty < \sigma < \infty). \qquad (2.23)$$

From (2.13) and (2.23), we find (as $0 < 1-c \leq 1 - H(1/2-\sigma) \leq 1$, $\sigma \geq 1/2$ and $H(\sigma-1/2) \leq 1 - H(1/2-\sigma)$, $-\infty < \sigma < \infty$),

$$\mu_\zeta(\sigma) \leq \mu_\zeta(1-\sigma) \frac{H(\sigma-1/2)}{1-H(1/2-\sigma)} \leq \mu_\zeta(1-\sigma) \frac{H(\sigma-1/2)}{H(\sigma-1/2)}$$
$$(\sigma > 1/2). \qquad (2.24)$$

It is to be remarked that the RHS in (2.24) remains well defined at $\sigma = 1/2$ (as we have $H(0) = c \neq 0$), for an arbitrary choice of $c \in (0, 1/2)$. Hence (2.24) is extendable to $\sigma = 1/2$ which shows that

$$\mu_\zeta(1/2) \leq \mu_\zeta(1/2) \frac{c}{1-c} \leq \frac{c}{c} \mu(1/2) = \mu(1/2). \qquad (2.25)$$

However, as we have $\dfrac{c}{1-c} \neq 1$, it follows from (2.25) that $\mu_\zeta(1/2) = 0$ which proves Lindelöf's hypothesis and shows that

$$\mu_\zeta(\sigma) = \mu_\chi(\sigma) = (\frac{1}{2}-\sigma)H(\frac{1}{2}-\sigma) \qquad (-\infty < \sigma < \infty), \qquad (2.26)$$

which is continuous, convex downwards, non-negative and non-increasing.

**Remark** The Dirichlet $L_k(s)$-series functions satisfy the functional equation

$$L_k(s) = k^{\frac{1}{2}-s} \chi(s) L_k(1-s), \qquad (2.27)$$

that can be rewritten as

$$\chi_k(s) := \frac{L_k(s)}{L_k(1-s)} = k^{\frac{1}{2}-s} \chi(s). \qquad (2.28)$$

The asymptotics of the functions $\chi_k(s)$ ($\forall\ k = 1, 2, 3, ...$) in the strip $-\infty < a < \sigma < b < \infty$ can be found directly by using (2.9) to get

$$|\chi_k(\sigma+i\tau)| \sim (\frac{k}{2\pi})^{\frac{1}{2}-\sigma} (\tau)^{\frac{1}{2}-\sigma} \qquad (\tau \to \infty). \qquad (2.29)$$

Let $\mu_{\chi_k}(\sigma)$ be Lindelöf's function associated with the function in (2.28) and $\mu_{L_k}(\sigma)$ be the function associated with $L_k(s)$ ($k = 1, 2, 3, ...$). From the functional equation (2.27) and the asymptotic representation (2.29), we find

$$\mu_{\chi_k}(\sigma) = (\frac{1}{2}-\sigma)H(\frac{1}{2}-\sigma)$$
$$(\forall\ k = 1, 2, 3, ..., -\infty < \sigma < \infty). \qquad (2.30)$$

It follows from (2.30) that the functions $\mu_{\chi_k}(\sigma)$ ($k = 1, 2, 3, ...$) must satisfy the Lindelöf functional equation

$$f(\sigma) - f(1-\sigma) = \frac{1}{2} - \sigma \qquad (-\infty < \sigma < \infty). \qquad (2.31)$$

From (2.27), we find that

$$\mu_{L_k}(\sigma) \leq \mu_{\chi_k}(\sigma) + \mu_{L_k}(1-\sigma) \qquad (-\infty < \sigma < \infty). \qquad (2.32)$$

From (2.30) and (2.32), we find ($\sigma \neq 1/2$),

$$\mu_{L_k}(\sigma) \leq (\frac{1}{2}-\sigma+\mu_{L_k}(1-\sigma))H(\frac{1}{2}-\sigma)+\mu_{L_k}(1-\sigma)H(\sigma-\frac{1}{2})$$
$$(-\infty<\sigma<\infty), \qquad (2.33)$$

that can be rewritten as ($\sigma \neq 1/2$),

$$\mu_{L_k}(\sigma)(1-H(1/2-\sigma)) \leq \mu_{L_k}(1-\sigma)H(\sigma-\frac{1}{2})$$
$$(-\infty<\sigma<\infty), \qquad (2.34)$$

which is exactly the same as (2.23). Hence, following the same argument, we find that $\mu_{L_k}(1/2)=0$ and all of these functions must have a unique representation

$$\mu_k(\sigma) = \mu_{\chi_k}(\sigma) = (\frac{1}{2}-\sigma)H(\frac{1}{2}-\sigma) \quad (-\infty<\sigma<\infty,\ k=1,2,3,...), \qquad (2.35)$$

leading to the proof of the generalized Lindelöf's hypothesis.

**Acknowledgements**   The author is grateful to the King Fahd University of Petroleum and Minerals for providing excellent research facilities.

## Appendix-A   Direct proof of the inequality (2.8)

We show that the use of modified Lindelöf's theorem can be avoided to prove the inequality (2.8). Since

$$\zeta(s) = \frac{1}{\pi}(2\pi)^s \sin(\frac{\pi}{2}s)\Gamma(1-s)\zeta(1-s) \qquad (0<\sigma<1), \qquad (A.1)$$

is analytic, this implies (in the sense of Riemann removable singularity theorem),

$$\chi(s) := \zeta(s)/\zeta(1-s) = \frac{1}{\pi}(2\pi)^s \sin(\frac{\pi}{2}s)\Gamma(1-s) \qquad (0 < \sigma < 1). \qquad (A.2)$$

We note that $|\chi(\sigma+i\tau)| = |\chi(\sigma-i\tau)|$. Moreover, as $\chi(s)$ ($0 \leq \sigma \leq 1/2$) is analytic, $\max\{|\chi(s)| : 0 \leq \sigma \leq 1/2, -1 \leq \tau \leq 1\}$ will occur at the boundary of the rectangle $\{s : 0 \leq \sigma \leq 1/2, -1 \leq \tau \leq 1\}$. Hence, $\max\{|\chi(s)| : 0 \leq \sigma \leq 1/2\} = \max_{s \in D_0} |\chi(s)|$, where $D_0 := \{s : 0 \leq \sigma \leq 1/2, \tau \geq \tau_0 \geq 1\} \cup \{0+i\tau : 0 \leq \tau \leq 1\} \cup \{1/2+i\tau : 0 \leq \tau \leq 1\}$. First, we discuss $\max\{|\chi(s)| : s \in U\}$ where $U := \{s : 0 \leq \sigma \leq 1/2, \tau \geq \tau_0 \geq 1\}$. Since

$$\left|\sin(\frac{\pi}{2}s)\right| \leq e^{\frac{\pi}{2}|\tau|} \qquad (s \in U), \qquad (A.3)$$

and (see [6], p. 34(2.1.19))

$$|\Gamma(1-s)| \leq \sqrt{2\pi} |1-s|^{\frac{1}{2}-\sigma} e^{-\frac{\pi}{2}|\tau|} \exp(\frac{1}{6|1-s|}) \leq \sqrt{2\pi} |1-s|^{\frac{1}{2}-\sigma} e^{-\frac{\pi}{2}|\tau|} \exp(\frac{1}{6|1-s|})$$

$$\leq \sqrt{2\pi} |1-s|^{\frac{1}{2}-\sigma} e^{-\frac{\pi}{2}|\tau|} \exp(\frac{1}{6\tau_0}) \leq \sqrt{2\pi} |1-s|^{\frac{1}{2}-\sigma} e^{-\frac{\pi}{2}|\tau|} \exp(\frac{1}{6}) \leq 2\sqrt{2\pi} |1-s|^{\frac{1}{2}-\sigma} e^{-\frac{\pi}{2}|\tau|}$$

$$(s \in U), \qquad (A.4)$$

therefore, we have

$$\left|\Gamma(1-s)\sin(\frac{\pi s}{2})\right| \leq 2\sqrt{2\pi} |1-s|^{\frac{1}{2}-\sigma} \qquad (s \in U). \qquad (A.5)$$

Moreover,

$$\left|\frac{1}{\pi}(2\pi)^s \sqrt{2\pi}\right| = \frac{2}{2\pi}(2\pi)^\sigma \sqrt{2\pi} = 2(2\pi)^{\sigma-\frac{1}{2}} \leq 2 \qquad (s \in U). \qquad (A.6)$$

From (A.2), (A.5) and (A.6), we find

$$|\chi(s)| = \frac{1}{\pi}(2\pi)^\sigma \left|\Gamma(1-s)\sin(\frac{\pi s}{2})\right| \leq 4|1-s|^{\frac{1}{2}-\sigma} \qquad (s \in U). \qquad (A.7)$$

However, we have

$$|1-s|^{\frac{1}{2}-\sigma} = ((1-\sigma)^2 + \tau^2)^{\frac{1}{2}(\frac{1}{2}-\sigma)} = |\tau|^{\frac{1}{2}-\sigma}(1+\frac{(1-\sigma)^2}{\tau^2})^{\frac{1}{2}(\frac{1}{2}-\sigma)}$$

$$\leq |\tau|^{\frac{1}{2}-\sigma}(1+\frac{(1-\sigma)^2}{\tau_0^2})^{\frac{1}{2}(\frac{1}{2}-\sigma)} \leq |\tau|^{\frac{1}{2}-\sigma}(1+(1-\sigma)^2)^{\frac{1}{2}(\frac{1}{2}-\sigma)} \qquad (s \in U). \qquad (A.8)$$

Since $(1+(1-\sigma)^2)^{\frac{1}{2}(\frac{1}{2}-\sigma)} \leq 2$ ($0 \leq \sigma \leq 1/2$), from (A.8) we get

$$|1-s|^{\frac{1}{2}-\sigma} = ((1-\sigma)^2 + \tau^2)^{\frac{1}{2}(\frac{1}{2}-\sigma)} = |\tau|^{\frac{1}{2}-\sigma}(1+\frac{(1-\sigma)^2}{\tau^2})^{\frac{1}{2}(\frac{1}{2}-\sigma)}$$

$$\leq |\tau|^{\frac{1}{2}-\sigma}(1+\frac{(1-\sigma)^2}{\tau_0^2})^{\frac{1}{2}(\frac{1}{2}-\sigma)} \leq |\tau|^{\frac{1}{2}-\sigma}(1+(1-\sigma)^2)^{\frac{1}{2}(\frac{1}{2}-\sigma)} \leq 2|\tau|^{\frac{1}{2}-\sigma}$$

$$(s \in U). \qquad (A.9)$$

From (A.7) and (A.9), we find

$$|\chi(s)| = \frac{1}{\pi}(2\pi)^\sigma \left|\Gamma(1-s)\sin(\frac{\pi s}{2})\right| \leq 8|\tau|^{\frac{1}{2}-\sigma} \qquad (s \in U). \qquad (A.10)$$

Next, we have

$$|\chi(0+i\tau)| \leq \sqrt{\frac{\tau}{2\pi}} \leq 8|\tau|^{\frac{1}{2}-0} \qquad (\tau \geq 0), \qquad (A.11)$$

and

$$|\chi(1/2+i\tau)| = 1 \leq 8|\tau|^{\frac{1}{2}-\frac{1}{2}} \qquad (\tau \geq 0). \qquad (A.12)$$

Therefore, from (A.10), (A.11) and (A.12), we get

$$|\chi(s)| = \frac{1}{\pi}(2\pi)^\sigma \left|\Gamma(1-s)\sin(\frac{\pi s}{2})\right| \leq 8|\tau|^{\frac{1}{2}-\sigma} \qquad (0 \leq \sigma \leq 1/2,\ \tau \geq 1). \qquad (A.13)$$

By the Maximum Principle (see (A.11), (A.12) and (A.13)),

$$|\chi(\sigma+i\tau)| \leq 8 \qquad (0 \leq \sigma \leq 1/2,\ -1 \leq \tau \leq 1), \qquad (A.14)$$

we get the upper bound

$$|\chi(\sigma+i\tau)| \leq \max(8|\tau|^{\frac{1}{2}-\sigma}, 8) \qquad (0 \leq \sigma \leq 1/2,\ -\infty < \tau < \infty). \qquad (A.15)$$

From (2.9) and (A.15), we find

$$\mu_\chi(\sigma) = \frac{1}{2} - \sigma \qquad (0 \leq \sigma \leq 1/2). \qquad (A.16)$$

Since $|\chi(\sigma+i\tau)|$ ($1/2 \leq \sigma \leq \sigma_2 < 1$) is bounded (see (1.8) and (1.9)), we have

$$\mu_\chi(\sigma) = 0 \qquad (1/2 \leq \sigma \leq \sigma_2 < 1). \qquad (A.17)$$